\documentclass{article}

\def \T {{\mathbf {T}}}

\def \Q {{\bf Q}}

\def \P {{\bf P}}

\usepackage[T1]{fontenc} 
\usepackage[cp1251]{inputenc} 
\usepackage[russian]{babel}

\begin{document}

  \rm \Large

\begin{center}   
{ \LARGE \bf
О спектральных задачах Колмогорова  и Рохлина \\ в классе  перемешивающих автоморфизмов
\vspace{5mm}

В.В. Рыжиков}
\end{center} 

\vspace{5mm}
 \it Let  $T$ be a staircase rank-one  construction with parameters $r_j \sim j^d$, $0<d<0.2$, then  its spectrum does not have the group property, and the product $T\otimes T$ has  homogeneous spectrum of  multiplicity 2. \rm

\vspace{5mm}
А. Н. Колмогоров еще давно высказывал предположение (см. например, \cite{Si}) о том, что максимальный спектральный тип эргодического автоморфизма всегда подчиняет свою свертку. Это свойство, представляющее собой
естественный континуальный аналог группового свойства спектра эргодического автоморфизма с дискретным спектром, было доказано Я. Г. Синаем для
специального класса автоморфизмов, удовлетворяющих условию А (см. \cite{Si}).
Оказалось, что в общем случае указанная гипотеза неверна. В настоящем
параграфе будут построены эргодические автоморфизмы, максимальные спектральные типы которых не обладают «групповым свойством» {} в указанном
выше смысле. (Цитата из обзора  \cite{KS})

\vspace{5mm}
Тема отсутствия группового  свойства спектра автоморфизма получила развитие: В.И. Оселедец построил слабо перемешивающий пример   \cite{O},  А.М. Степин установил типичность отсутствия  свойства Колмогорова  \cite{S}. Эти примеры в силу их определения не обладали свойством перемешивания.  Переход к перемешивающим примерам произошел благодаря другой задаче, 
известной как  проблема Рохлина о непростом однородном спектре (об истории и эволюции этой задачи см. \cite{An},  \cite{Da}).
 
В \cite{99} автором было  доказано  существование лестничной конструкции  $T$ такой, что произведение  $T\otimes T$  имеет спектр кратности 2.  (Автоморфизм пространства Лебега $(X,\mu)$ и  индуцированный им унитарный оператор в $L_2(X,\mu)$    обозначаются одинаково.)      В этом случае ее спектральная  $\sigma_T$ и  сверточный квадрат $\sigma_T^{\ast 2}$  взаимно сингулярны. С учетом того, что  Т. Адамс \cite{Ad} доказал свойство перемешивания подходящего класса лестничных конструкций, в \cite{99} были решены сразу две задачи  в классе перемешивающих автоморфизмов: доказано отсутствие  группового свойства  спектра  для некоторых лестничных конструкции, и реализован 
  однородный  спектр кратности 2. Однако,   в силу специфики доказательства, основанном на "бесконечно малых" слабых пределах,   явные примеры не были  указаны. Это относится и к работе \cite{Ag}, где использовалась методика  \cite{99} (приближение перемешивающих систем неперемешивающими с нужным спектральным свойством).  

 Конкретные автоморфизмы появились  в \cite{20}. 
Лестничная конструкция задается  последовательностью параметров $r_j$.    В \cite{20}  показано, что  нужными примерами служат 
все конструкции, для которых  $r_j \sim \log j$.  Давайте, расширим класс  примеров с соответствующими пояснениями. \rm

\vspace{3mm}
\bf Теорема 1. \it Если $T$ -- лестничная конструкция с параметрами $r_j \sim j^d$, $0<d<0.2$,  то ее спектр не обладает групповым свойством, а произведение  $T\otimes T $ имеет  однородный  спектр кратности 2.   \rm

\vspace{3mm}
Существование перемешивающих автоморфизмов с однородным спектром кратности $n>2$ установлено   в \cite{T}. Конструктивного решения пока нет.

\section{Лестничная конструкция}
Пусть $r_j\to\infty$ и для всех $j$, начиная с некоторого,  
$\bar s_j=(1,2,\dots, r_j-2,r_j-1,0).$  

\it Определение лестничной конструкции.  \rm Пусть на шаге $j$  дана система   непересекающихся полуинтервалов 
$$E_j, TE_j, T^2E_j,\dots, T^{h_j-1}E_j,$$
причем на полуинтервалах $E_j, TE_j, \dots, T^{h_j-2}E_j$
пребразование $T$ является параллельным переносом. Такой набор   полуинтервалов  называется башней этапа $j$, их объединение обозначается через $X_j$ и также называется башней.

Представим   $E_j$ как дизъюнктное объединение  $r_j$ полуинтервалов 
$$E_j^1,E_j^2E_j^3,\dots E_j^{r_j}$$ одинаковой длины.  
Для  каждого $i=1,2,\dots, r_j$ определим  колонну $X_{i,j}$ как объединение интервалов  
$$E_j^i, TE_j^i ,T^2 E_j^i,\dots, T^{h_j-1}E_j^i.$$
К каждой  колонне $X_{i,j}$ добавим  $s_j(i)$  непересекающихся полуинтервалов  той же меры, что у $E_j^i$, получая набор  
$$E_j^i, TE_j^i, T^2 E_j^i,\dots, T^{h_j-1}E_j^i, T^{h_j}E_j^i, T^{h_j+1}E_j^i, \dots, T^{h_j+s_j(i)-1}E_j^i$$
(все эти множества  не пересекаются).
Обозначив $E_{j+1}= E^1_j$, для   $i<r_j$ положим 
$$T^{h_j+s_j(i)}E_j^i = E_j^{i+1}.$$
 Набор надстроеных колонн с этого момента  рассматривается как   башня  этапа $j+1$,  состоящая из полуинтервалов  
$$E_{j+1}, TE_{j+1}, T^2 E_{j+1},\dots, T^{h_{j+1}-1}E_{j+1},$$
где  
 $$ h_{j+1} =\sum_{i=1}^{r_j-1}(h_j+s_j(i)).$$

Частичное определение преобразования $T$ на этапе $j$ сохраняется на всех следующих этапах. На пространстве  $X=\cup_j X_j$ тем самым определено  обратимое преобразование $T:X\to X$, сохраняющее  стандартную меру Лебега на $X$.

Лестничная конструкция  эргодичена,   имеет простой спектр.  Известно, что индикаторы  интервалов $E_j$ являются циклическими векторами для оператора $T$.

Фиксируем индикатор $f$  некоторого этажа. 
Докажем, что для выбранной лестничной  конструкции $T$  для всех $r>0$  векторы $\T^rf\otimes f + f\otimes \T^rf$
принадлежат циклическому пространству $C_{f\otimes f}$  оператора $T\otimes T$.
Это означает, что симметрическая степень $T\odot T$ имеет простой  спектр, откуда вытекает отсутствие группового свойства спектра оператора $T$ и    двукратность спектра
произведения $T\otimes T$.
На самом деле нам достаточно показать, что для расстояния
$\rho$ между вектором в $L_2$ и подпространством выполнено   
$$\rho(T^rf\otimes f + f\otimes T^rf\,,\, C_{f\otimes f})\to 0.$$
Отсюда, как известно, вытекает, что $\T^rf\otimes f + f\otimes \T^rf$
принадлежат  $C_{f\otimes f}$.

Обозначим  $$ Q_r = \frac 1 r\sum_{i=0}^{r-1} T^{-i}, \ \ \ \ \Q_r=Q_r\otimes Q_r.$$ 
Имеет место  равенство
$$T^rf\otimes f + f\otimes T^rf
=r^2Q_rf\otimes Q_rf \ \ \  + \ \ \ (r-2)^2
 Q_{r-2}Tf\otimes Q_{r-2}Tf -\eqno (1)$$
$$-(r-1)^2
 Q_{r-1}f\otimes Q_{r-1}f \  
- \ (r-1)^2
 Q_{r-1}Tf\otimes Q_{r-1}Tf.$$

Предполагая  последовательность  $r_j\to\infty$ монотонной, положим 
$$ J_r=\{j:\, r_j=r+1,  \ \  j<j_r-r\},$$
где $j_r=\max\{j:\, r_j=r+1,\}$.

Мы найдем такое $D>0$, что при   $|J_r|> r^D$ будет выполнено 
$$\left\|\Q_r(f\otimes f)- \P_r(f\otimes f)\right\|^2 < 
\,  \frac 2 {|J_r|}=o\left(\frac 1{r^4}\right),  \eqno (2)$$
где
$$ \P_r(f\otimes f) = \frac 1 {|J_r|}\sum_{j\in J_r} T^{h_j}f\otimes T^{h_j}f.$$ 
Если  при    $r\to\infty $   выполнено    
норма разности $\P_r(f\otimes f)-\Q_r(f\otimes f)$ является $o\left(\frac 1{r^2}\right),$ то с учетом  
равенства $(1)$ получим, что векторы $T^rf\otimes {\bf 1} + {\bf 1}\otimes T^r$ лежат в циклическом пространстве вектора $f\otimes f$, что нам и надо. 

Рассмотрим  равенство
$$ {|J_r|^2}\|\Q_r(f\otimes f)-\P_r(f\otimes f)\|^2 =\ {|J_r|^2}(\Q_r(f\otimes f),\Q_r(f\otimes f))-  $$
$$ - 2{|J_r|}\left(\Q_r(f\otimes f)\,,\,\sum_{j\in J_r} T^{h_j}f\otimes T^{h_j}f\right)\ +\
\left(\sum_{j\in J_r} T^{h_j}f\otimes T^{h_j}f
\,,\,\sum_{k\in J_r} T^{h_k}f\otimes T^{h_k}f\right). \eqno (3)$$
Если значения скалярных произведений $\left(\Q_r(f\otimes f)\,,\, T^{h_j}f\otimes T^{h_j}f\right)$ достаточно близки к величине  $(\Q_r(f\otimes f),\Q_r(f\otimes f))$ и к ней же близки значения  
$\left( T^{h_j}f\otimes T^{h_j}f\,,\,T^{h_k}f\otimes T^{h_k}f\right)$ при $j\neq k$, 
то будет выполняться нужное нам неравенство $(2)$.  Теперь нужно пояснить, что в нашем случае
означает "достаточно близки". 

\vspace{3mm}
\bf Лемма. \it Пусть  $j\in J_r$.
Тогда $$|(T^{h_j}f, Q_{r}f) - (Q_{r}f,Q_{r}f)|< \frac {2r}{h_j}+ 2r^{-r}=o(r^{-1000}). $$\rm

\vspace{3mm}
При $j,\in J_r$   из леммы вытекает   
 $$|\left(\Q_r(f\otimes f)\,,\, T^{h_j}f\otimes T^{h_j}f\right)\,  -\, (\Q_r(f\otimes f),\Q_r(f\otimes f))|=o(r^{-100}).$$ 
Для некоторой константы $C$ для  $j, k=j+p\in J_r$, $p\geq 0$,  также имеем
$$|\left(T^{h_j}f\otimes T^{h_j}f\,,\, T^{h_{j+p}}f\otimes T^{h_{j+p}}f\right)\,  -\, (\Q_r(f\otimes f),\Q_r(f\otimes f))|< \frac {Ch_j}{h_{j+p}}\leq Cr^{-p}.$$

В  $(3)$ основной вклад  в правую часть равенства  дают разности 
$$\left(T^{h_j}f\otimes T^{h_j}f
\,,\,T^{h_{j+p}}f\otimes T^{h_{j+p}}f\right) -\left(\Q_r(f\otimes f)\,,\, T^{h_k}f\otimes T^{h_k}f\right)
$$
при $p=0$.
 Таким образом, при   $r_j\sim j^d$,  $0<d<\frac 1 5$,  имеем 
$|J_r|\sim  r^{\frac {1-d} d}$ и $$  \frac 1 {|J_r|}=o\left(\frac 1{r^4}\right).$$ Это приводит, как мы пояснили, к простому спектру произведения $T\odot T$ для соответствующей лестничной конструкции $T$.

В случае пространства с сигма-конечной мерой "перемешивающие" {} автоморфизмы (т.е. $T^n\to_w 0$)  без группового свойства спектра явно предъявляются  среди сидоновских конструкций \cite{24}.

\end{document}